\documentclass[12pt]{article}

\usepackage{amsfonts}
\usepackage{amssymb}

\usepackage{amsmath,amssymb}
\usepackage{epsfig,afterpage}
\usepackage{graphicx}
\usepackage{amsmath}
\usepackage{amsthm}
\usepackage{amscd}
\usepackage{amsbsy}
\usepackage{epsfig}
\usepackage[dvips]{psfrag}
\usepackage{makeidx}

\usepackage[latin1]{inputenc}

\newtheorem{theorem}[subsection]{Theorem}

\newtheorem{lemma}[subsection]{Lemma}

\newtheorem{proposition}[subsection]{Proposition}
\newcommand{\N}{\mathbb{N}}

\newcommand{\R}{\mathbb{R}}

{ \small \normalsize}

\newfont{\bms}{msbm10 scaled 1200}
\newfont{\bmt}{msbm10 scaled 2500}

\newfont{\bmtt}{msbm10 scaled 1700}

\title{Analytic varieties as limit periodic sets}
\author{
        Andr\'e Belotto \\
        Universit\'e de Haute-Alsace\\
        andre-ricardo.belotto-da-silva@uha.fr
}
\date{07/2011}

\begin{document}
\maketitle

\section*{Abstract}
Let $f(x,y) \not\equiv 0$ be a real-analytic planar function. We show that, for almost every $R>0$ there exists an analytic $1$-parameter family of vector fields $X_{\lambda}$ which has $\{f(x,y)=0\} \cap \overline{B_R\left((0,0)\right)}$ as a limit periodic set. Furthermore, we show that if $f(x,y)$ is polynomial, then there exists a polynomial family with these properties.
\section{Introduction} \label{sec:in}
Hilbert, in $1900$, proposed to find an upper bound for the number of limit cycles depending on the degree of the planar real-polynomial vector field, and to analyze the `distribution' of the limit cycles on the plane. This problem is known as the $16^o$ Hilbert problem. Ilyashenko $\cite{IL}$ and Écalle $\cite{Ec}$ proved, independently, that any polynomial vector field has a finite number of limit cycles. But, it is still an open question whether there exists an upper bound for the number of limit cycles depending on the degree of the vector field.\\
\\
Roussarie proposes a program to show the existence of this upper bound (see e.g. $\cite{Rou}$). This program is based on the notion of limit periodic sets for analytic families of planar vector fields, which was first defined by Françoise-Pugh in $\cite{Fr}$. We recall that an \textit{analytic family of planar vector fields} is given by an open set $U \subset \R^2$, a $n$-dimensional analytic manifold $\Lambda$ (called the parameter space) and, for each $\lambda \in \Lambda$, a vector field $X_{\lambda}= a(x,y,\lambda) \frac{\partial}{\partial x} + b(x,y,\lambda) \frac{\partial}{\partial y} $ defined in $U$ such that the functions $a,b : U \times \Lambda \rightarrow \R$ are analytic. If $\Lambda$ is an open set of $\R^n$ and $a,b$ are polynomials, then we say that $X_{\lambda}$ is a polynomial family of planar vector fields.\\
\\
\textbf{Definition:}($\cite{Fr}$) A limit periodic set for $X_{\lambda}$ is a compact non-empty subset $\Gamma \subset U$, such that there
    exists a sequence $(\lambda_n)_{n \in \mathbb{N}} \rightarrow \lambda_0$ in the
    parameter space, and a
    limit cycle $\gamma_{n} \subset U$ of $X_{\lambda_n}$ such that $\gamma_{n} \rightarrow \Gamma$ (in the Haussdorf metric) as $n \rightarrow \infty$.\\
    \\
The study of limit periodic sets is a classical subject in bifurcation theory. Some usual non-trivial examples are: the homoclinic loop; the weak focus; the heteroclinic graph. In relation with Roussarie's program, Roussarie and Panazzolo give in $\cite{RP}$ an explicit example of a polynomial family with a segment as limit periodic set. A natural question which arises is the following:\\
\\
\textbf{Problem:} which kind of subset $\Gamma \subset U$ can occur as limit periodic set of an analytic (or polynomial) family of vector fields?\\
\\ 
In this paper we will go further on the study of the structure of limit periodic sets. Our main result provides many new examples of degenerated limit periodic sets. Given an analytic function $f: U \subset \R^2 \rightarrow \R$ and a positive scalar $R\in \R$, we consider the following hypotheses:
\begin{itemize}
\item[i] The closure of the ball $B_R=\{p \in \R^2; ||p|| < R \} $ is contained in $U$;
\item[ii] The zero set of $f$, $Z(f)=f^{-1}(0)$, and the closed ball $\bar{B}_R$ has non-empty intersection. We note by $\Gamma=Z(f) \cap \bar{B}_R$ this intersection;
\item[iii] The function $\nabla f(x,y)$ restricted to $\Gamma$, is non zero outside a finite set of points;
\item[iv] The set $Z(f)$ intersects the boundary of $\bar{B}_R$ in a finite set of points, and the function $<\nabla f(p)^{\bot},p>$ is non-zero at these points.
\end{itemize}
%
The main result of the paper is:
\begin{theorem}
Let $f: U \subset \R^2 \rightarrow \R$ be an analytic function and $R\in \R^+$ such that hypotheses $[i]-[iv]$ are satisfied. Then, there exists an analytic family $X_{\lambda}$, with parameter space $\Lambda= \mathbb{R}$, such that each connected component of $\Gamma$ is a limit periodic set at $\lambda=0$. Furthermore, if $f(x,y)$ is polynomial of degree $M$, then we can find $X_{\lambda}$ polynomial of degree at most $2(M-1)(M-2)+7$.
\label{th:main}
\end{theorem}
With weaker hypotheses, we obtain the following result:
\begin{theorem}
Let $f: U \subset \R^2 \rightarrow \R$ be a non-zero analytic function and $R\in \R^+$ such that the hypotheses $[i]-[iii]$ are satisfied. Then there exists an analytic family $X_{\lambda}$, with parameter space $\Lambda= \mathbb{R}$, such that $\Gamma$ is a union of limit periodic sets at $\lambda=0$. Furthermore, if $f(x,y)$ is  polynomial of degree $M$, then we can find $X_{\lambda}$  polynomial of degree at most $2(M-1)(M-2)+7$.
\label{th:1}
\end{theorem}
\textbf{Remark:}
The bound for the degree of $X_{\lambda}$ obtained in theorem $\ref{th:1}$ is not optimal, as we obtain it as a corollary of the more elaborate proof of theorem $\ref{th:main}$. Indeed, we could prove the same theorem with only hypotheses $[i]$ and $[ii]$, and show that, if $f(x,y)$ is polynomial of degree $M$, we can construct $X_{\lambda}$  polynomial of degree exactly $4M-1$. Through the paper we will briefly indicate the necessary changes to improve this result.\\
\\
In a forthcoming paper we will study the topology of arbitrarily degenerated limit periodic sets using a desingularization processes.
\subsection{Examples}
\begin{figure}
  \centering
  \includegraphics[scale=0.5]{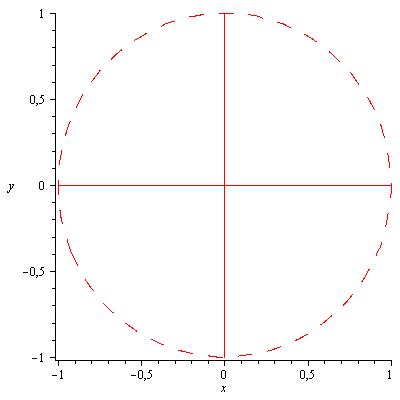}
\includegraphics[scale=0.5]{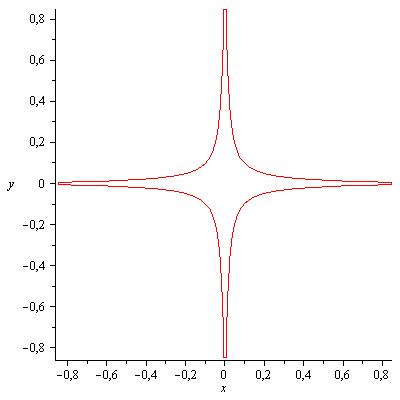}
  \caption{\small{At left, $\{f(x,y)=0\}$ and $\overline{B_1((0,0))}$. At right, the limit cycle for $\lambda=\frac{1}{100}$.}}
  \label{fig:ex1a}
 \end{figure}
\textbf{First example:} Let us consider $f=xy$, so $Z(f)=\{x=0\} \cup \{y=0\}$ is a `cross' (see figure $\ref{fig:ex1a}$ at left). It is clear that, for any $R>0$ all conditions of theorem $\ref{th:main}$ are satisfied, in particular, we take $R=1$. Now consider the vector field:
\begin{equation}
X_{\lambda}= \left\{
\begin{array}{cc}
\dot{x}= & {\lambda}y+y{x}^{2}+ \left( {\lambda} \left( {x}^{2}
+{y}^{2}-1 \right) +{x}^{2}{y}^{2} \right)  \left( {\lambda}x+x{y}
^{2} \right) 
\\
\dot{y}= & -{\lambda}x-x{y}^{2}+ \left( {\lambda} \left( {x}^{2
}+{y}^{2}-1 \right) +{x}^{2}{y}^{2} \right)  \left( {\lambda}y+y{x
}^{2} \right) 
\end{array}
\right.
\label{eq:excruz}
\end{equation}
For each $\lambda>0$, $X_{\lambda}$ has an unique limit cycle given by:
$$
\gamma_{\lambda}=\{(xy)^2+\lambda(x^2+y^2-1)=0\}
$$
which converges to $Z(f) \cap \bar{B}_1$ as $\lambda$ goes to zero. Figure $\ref{fig:ex1a}$ at right shows a limit cycle of this field.\\
\\
 \begin{figure}
  \centering
  \includegraphics[scale=0.5]{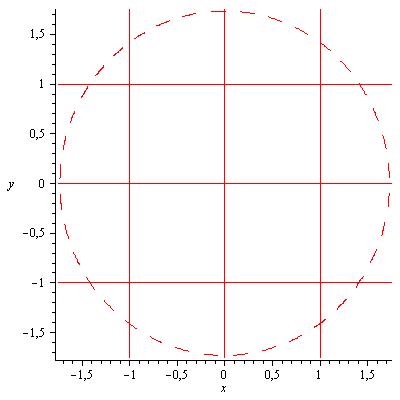}
\includegraphics[scale=0.5]{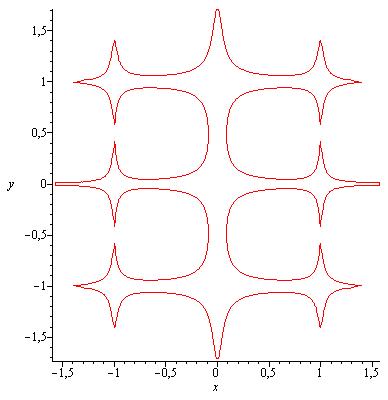}
  \caption{\small{At left, $Z(f)$ and $\bar{B}_R$. At right, the limit cycle for $\lambda=\frac{1}{10000}$.}}
  \label{fig:triv}
 \end{figure}
\textbf{Second example:} Let us consider $f=xy(x+1)(x-1)(y+1)(y-1)$ and $R=\sqrt{3}$. Then $Z(f) \cap \bar{B}_R$ is represented by the figure $\ref{fig:triv}$ at left. We now construct the vector field (as given in section $\ref{sec:2}$): 
$$ 
X_{\lambda} =  X_{H_{\lambda}}+ H_{\lambda}(x,y) \nabla H_{\lambda}(x,y)
$$
where $X_{H_{\lambda}}$ is the Hamiltonian of $H_{\lambda}(x,y)$ and:
$$
H_{\lambda}(x,y) = f(x,y)^2+ \lambda(x^2+y^2-R^2) \prod_{e \in E}{((x-x_e)^2+(y-y_e)^2-\lambda^2)}
$$
where $E=\{1,2,3,4\}$, $p_1=(1,\frac{1}{2})$, $p_2=(1,-\frac{1}{2})$, $p_3=(-1,\frac{1}{2})$ and $p_4=(-1,-\frac{1}{2})$.\\
\\
For each $\lambda>0$, $X_{\lambda}$ has an unique limit cycle given by:
$$
\begin{array}{c}
\gamma_{\lambda} = \{f(x,y)^2+\lambda(x^2+y^2-1)((x+1)^2+(y+\frac{1}{2})^2-\lambda^2)((x+1)^2+(y-\frac{1}{2})^2-\lambda^2)\\
((x-1)^2+(y+\frac{1}{2})^2-\lambda^2)((x-1)^2+(y-\frac{1}{2})^2-\lambda^2)=0\}
\end{array}
$$
which converges to $Z(f) \cap \bar{B}_{\sqrt{3}}$ as $\lambda$ goes to zero. Figure $\ref{fig:triv}$ at right, shows a limit cycle of this field. Notice that it is connected.\\ 
\\ 
 \begin{figure}
  \centering
\includegraphics[scale=0.5]{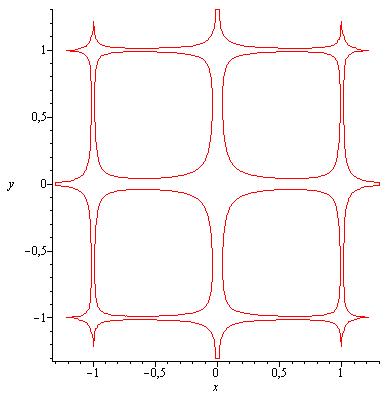}
  \caption{\small{The limit cycles for $\lambda=\frac{1}{10000}$ of the third example. Remark the contrast with figure $\ref{fig:triv}$ at right.}}
  \label{fig:triv2}
 \end{figure}
Let us now illustrate what is obtained applying the weaker result of theorem $\ref{th:1}$.\\
\\
\textbf{Third example:} Again let us consider $f=xy(x+1)(x-1)(y+1)(y-1)$ and $R=\sqrt{3}$. Take the vector field:
$$ 
X_{\lambda} =  X_{H_{\lambda}}+ H_{\lambda}(x,y) \nabla H_{\lambda}(x,y)
$$
where:
$$
H_{\lambda}(x,y) = f(x,y)^2+ \lambda(x^2+y^2-R^2)
$$
For each $\lambda > 0$ small enough, $X_{\lambda}$ has five limit cycles (see figure $\ref{fig:triv2}$) given by the five connected components of the set:
$$
\{f(x,y)^2+\lambda(x^2+y^2-1)=0\}
$$
which converges to $Z(f) \cap \bar{B}_{\sqrt{3}}$ as $\lambda$ goes to zero. In general, the polynomial families obtained through theorem $\ref{th:1}$ are of much smaller degree. However we can only guarantee the $\Gamma$ is a union of a finite number of limit periodic sets (which is five at this example).
\subsection{Sketch of the proof}
The basic idea to obtain the family $X_{\lambda}$ is to find a function $H(x,y,\lambda)$ with the following properties for $\lambda>0$:
\begin{itemize}
\item $G_{\lambda} = \{(x,y); H(x,y,\lambda)=0\}$ is compact;
\item The sets $G_{\lambda}$ converges (on the Haussdorf topology) to $\Gamma$ as $\lambda$ goes to zero;
\item There exists a sequence $(\lambda_n)_{n \in \N} \rightarrow 0$, such that $G_{\lambda_n}$ is regular 
; 
\item If $\Gamma$ is the union of $N$ connected components, then $G_{\lambda}$, for $\lambda>0$ sufficiently small, is the union of $N$ connected components, each one of them converging (on the Hausdorff topology) to a different connected component of $\Gamma$. 
\end{itemize}
Once we get this function, we will explicit construct the analytic-family $X_{\lambda}$ as a perturbation of its Hamiltonian. All this construction is done in section $\ref{sec:1}$.\\
\\
The section $\ref{sec:2}$ is devoted to the proof of the first three properties, and section $\ref{sec:3}$ to the proof of the fourth property. The fourth property is more complicated: the basic idea is to get a local description of $G_{\lambda}$ for $\lambda$ small enough and proceed with a topological global argument. 
%
%
%
\section{The explicit family} \label{sec:1}
At this section, assuming two technical propositions, we prove theorems $\ref{th:main}$ and $\ref{th:1}$ by explicitly constructing the analytic family $X_{\lambda}$. The proof of these technical propositions will be given in the remaining sections of the paper.\\
\\
From now on, we will always consider $f(x,y) \not\equiv 0$ and $R \in \R^{+}$ satisfying hypotheses $[i]-[iii]$. By the Newton-Puiseux theorem (see e.g. $\cite{Wall}$, chap 2) and the compactness of $\bar{B}_R$, it is clear that $B_R \setminus Z(f)$ is a finite union of connected open sets $C_j$, where $j \in J$ for some index set $J \subset \N$. By the same reasoning, it is easy to see that the boundary of each $C_j$ is a finite union of points and regular $1$-manifolds. Call $(C_j)_{j \in J}$ the \textit{first support of $(f,R)$}.\\
\\
A $C_j$ is \textit{exterior} if $\partial C_j \cap \partial \overline{B_R(p_0)}$ contains a $1$-manifold. It will be called \textit{interior} otherwise. Take the union: $\widetilde{C}=\cup_{k\in K} C_k$, where $K=\{k\in J; C_j $ is exterior $\}$. Defining $I=(J \setminus K )\cup \{0\}$ and indexing $\widetilde{C}$ by $C_0$, we define the \textit{support of $(f,R)$} as $(C_i)_{i \in I}$.\\
\\
An \textit{adjacency graph of $(f,R)$} is the simple graph $G=(I,F)$ where, $\{i,j\} \in F$ if $\partial C_i \cap \partial C_j$ contains a $1$-manifold. A spanning tree $T=(I,E)$ of $G$ 
will be called an \textit{adjacency tree of $(f,R)$} (for basic concepts of graph theory see e.g. $\cite{Graph}$).\\
\\
\textbf{Lemma:}
\textit{There exists an adjacency tree for $(f,R)$.}
\begin{proof}
It is a simple exercise to show that a graph is connected if, and only if, it admits a spanning tree. Take $I_0=\{0\}$, $I_1=\{ i \in I; \partial C_i \cap \partial C_0$ contains a $1-$manifold $\}$ and so on. As the number of components $C_i$ is finite, this process ends and we have that $I_{\infty} = \cup_{j \in \N} I_j$ is a connected subcomponent of $I$. Now, by the Newton-Puisseux theorem and maximality of $I_{\infty}$, the boundary of $\cup_{i \in I_{\infty}}\bar{C}_i$ must be $\partial \bar{B}_R$. So, $I=I_{\infty}$.
\end{proof}
Fix an adjacency tree $T=(I,E)$ of $(f,R)$. Now, for each $\{e_1,e_2\} \in E$, we chose a point $p_e=(x_e,y_e)$ such that:
\begin{itemize}
\item $p_e$ is on the $1$-dimensional boundary of $C_{e_1}$ and $C_{e_2}$;
\item $\nabla f(p_e) \neq (0,0)$. 
\end{itemize}
We will denote by $\Gamma_E$ the set of all $p_e$, for $e \in E$.\\
\\
\textbf{Remark:}
For each $e \in E$, the existence of $p_e$ is guaranteed by the hypotheses $[iii]$.\\
\\
Take now:
$$
h(x,y,\lambda) = f(x,y)^2+\lambda(x^2+y^2-R^2) \prod_{e \in E}{((x-x_e)^2+(y-y_e)^2-\lambda^2)}
$$
and $H(x,y,\lambda,\alpha) = h(x,y,\lambda) - \alpha \lambda^4$.\\
\\
\textbf{Remark:}
If $E=\emptyset$, then we convention that $\prod_{e \in E}{((x-x_e)^2+(y-y_e)^2-\lambda^2)}=1$.\\
\\
For a fixed $\alpha_0 $ we write $H_{\lambda, \alpha_0}(x,y) = H_{\alpha_0}(x,y,\lambda)=H(x,y,\lambda,\alpha_0)$ and $h_{\lambda}(x,y)=h(x,y,\lambda)$, and consider the $1$-parameter family of vector fields:
$$
X_{\lambda}=X_{\lambda,\alpha_0}= X_{H_{\lambda,\alpha_0}} + H_{\lambda,\alpha_0}(x,y) \nabla H_{\lambda,\alpha_0}(x,y)
$$
where $X_{H_{\lambda,\alpha_0}}$ stands for the Hamiltonian of $H_{\lambda,\alpha_0}(x,y)$.\\
\\
\textbf{Remark:}
If $f(x,y)$ is a polynomial, then $X_{\lambda}$ is polynomial. More than this, if $M$ is the degree of $f$, then, by Harnack's inequality (see e.g. $\cite{Har}$), the degree of $h_{\lambda}(x,y) $ is bounded by $ max\{2M, (M-1)(M-2)+2\}$. This imply that the degree of $X_{\lambda}$ is bounded by $2(M-1)(M-2)+7$.\\
\\
Denote by $G_{\lambda,\alpha}$ the set $\{(x,y); H_{\lambda,\alpha}(x,y)=0\}$. If there is no risk of confusion about the $\alpha$ fixed, simply denote $G_{\lambda,\alpha}$ by $G_{\lambda}$. The set $G_{\lambda,\alpha}$ is \textit{regular} if $\nabla_{x,y}(H_{\lambda,\alpha})(p) \neq (0,0)$ for all $p \in G_{\lambda,\alpha}$. 
\begin{proposition}
Let $f: U \subset \R^2 \rightarrow \R$ a non-zero analytic function and $R \in \R^{+}$ satisfies hypotheses $[i]-[iii]$. For almost all $\alpha \in [0,1]$ fixed and $\lambda>0$, $H_{\lambda,\alpha}(x,y)$ is such that:
\begin{itemize}
\item $G_{\lambda}$ is compact;
\item The set $G_{\lambda}$ converges (on the Hausdorff topology) to $\Gamma$ when $\lambda$ goes to zero;
\item There exists a sequence $(\lambda_n)_{n \in \N} \rightarrow 0$, such that $G_{\lambda_n}$ is regular. 
\end{itemize}
\label{prop:123}
\end{proposition}
\begin{proposition}
Keeping the same hypotheses of the preceding Proposition, if $f$ and $R$ further satisfy hypothesis $[iv]$, and $\Gamma$ is the union of $N$ connected components, then $G_{\lambda}$, for $\lambda>0$ sufficiently small, is the union of $N$ connected components, each one of them converging (on the Hausdorff topology) to a different connected component of $\Gamma$. 
\label{prop:4}
\end{proposition}
\textbf{Remark:}
If we take: 
$$
H^{\ast}(x,y,\lambda,\alpha) = f(x,y)^2+\lambda(x^2+y^2-R^2) - \alpha \lambda^2
$$
We claim that:
\begin{itemize}
\item Proposition $\ref{prop:4}$ would not be true (see third example of the introduction);
\item The proposition $\ref{prop:123}$ could be proven with only hypotheses $[i]$ and $[ii]$;
\item If $f(x,y)$ is a  polynomial of degree $M$, then the degree of $H^{\ast}(x,y,\lambda,\alpha)$ would be $2M$. This imply that $X_{\lambda}$ would have degree exactly $4M-1$.
\end{itemize}
As the main goal of the paper is to prove theorem $\ref{th:main}$, we leave the proof of this claim to the reader, since it is a simple adaptation of the proofs contained in section $\ref{sec:2}$.\\
\\
Assuming that the propositions are true, we can easily prove the theorems. Indeed:
\begin{proposition}
For $\lambda \neq 0$ and a fixed $\alpha \in [0,1]$:
\begin{itemize}
\item All periodic orbits of $X_{\lambda}=X_{\lambda,\alpha}$ are contained in $G_{\lambda}$;
\item If $G_{\lambda}$ is compact and regular then $G_{\lambda}$ is a union of limit cycles.
\end{itemize}
\label{prop:fXl}
\end{proposition}
\begin{proof}
Start noticing that:
$$
X_{\lambda} (H_{\lambda,\alpha}(x,y)) = H_{\lambda,\alpha}(x,y) || \nabla H_{\lambda,\alpha}(x,y)||^2
$$
So, if $p \notin G_{\lambda}$, the solution $\gamma(t)$ passing through it is such that:
\begin{itemize}
\item Either $H_{\lambda,\alpha}(\gamma(t))$ is a strictly increasing or decreasing function;
\item Or $\nabla H_{\lambda,\alpha}(p) = 0$ and $p$ is a singularity of $X_{\lambda}$.
\end{itemize}
In any case, $p$ can not belong to a periodic orbit.\\
\\
Now, if $G_{\lambda}$ is compact, as it it is invariant by $X_{\lambda}$, it must be the union of periodic orbits, singular points and orbits ending in singular points. Since $G_{\lambda}$ is regular, there are no singular points of $X_{\lambda}$ in $G_{\lambda}$, thus it is a union of periodic orbits. As there are no other periodic orbits and $H_{\alpha,\lambda}(x,y)$ is analytic, each such periodic orbit is necessarily a limit cycle.
\end{proof}
Now, clearly:
\begin{itemize}
\item Theorem $\ref{th:1}$ is a corollary of  propositions $\ref{prop:123}$ and $\ref{prop:fXl}$;
\item The main Theorem  $\ref{th:main}$ is a corollary of propositions $\ref{prop:123}$, $\ref{prop:4}$ and $\ref{prop:fXl}$.
\end{itemize}
\section{Proof of Proposition $\ref{prop:123}$} \label{sec:2}
The proof consists in three small lemmas:
\begin{lemma}
For $\lambda_0 > 0$ and $\alpha_0 \in [0,1]$, every level curve of $H_{\lambda_0,\alpha_0}(x,y)$ is compact.
\label{lem:compactLC}
\end{lemma}
\begin{proof}
Take $r \in \R$. Fixing $\lambda_0 \neq 0$ and $\alpha_0 \in [0,1]$. As $|E| < \omega$ and $f(x,y)^2 \geq 0$, and
$$
h(x,y,\lambda_0) = f(x,y)^2  + \lambda_0( \sum^{|E|+1}_{i=0} x^{2i} y^{2(|E|+1-i)} + \widetilde{h}(x,y,\lambda))
$$
where $\widetilde{h}(x,y,\lambda)$ is a  polynomial of degree smaller then $2|E|+2$, it is clear that there exists $\tau \in \R$ such that $h(x,y,\lambda_0)> \alpha_0 \lambda_0^4 + r$ for all $(x,y) \notin B_{\tau}$. So, $\{H_{\lambda_0,\alpha_0}(x,y)=r\}$ must be contained in $B_{\tau}$ and is limited. As $H$ is a continuous function we are done.
\end{proof}
\begin{lemma}
For $\lambda>0$, the set $G_{\lambda}$ converges (on the Haussdorf topology) to $\Gamma$ as $\lambda$ goes to zero.
\label{lem:convergencia}
\end{lemma}
\begin{proof}
We will proceed in four steps:\\
\\
Step [I]: As $\lambda$ goes to $0$, $G_{\lambda}$ must converge to a compact set $\gamma$. Indeed:
$$
H_{\alpha}(x,y,\lambda) = f(x,y)^2  + \lambda_0( \sum^{\|E|+1}_{i=0} x^{2i} y^{2(|E|+1-i)} + \widetilde{h}(x,y,\lambda) - \alpha_0 \lambda^3)
$$
where $\widetilde{h}(x,y,\lambda)$ is a  polynomial of degree in $x$ and $y$ smaller then $2|E|+2$. So, if $\lambda$ is small enough, $G_{\lambda}$ must be contained in $B_{2R} = \{p ; ||p||^2  < 2R \}$: all $G_{\lambda}$ are contained in a same compact set for small enough $\lambda$.\\
\\
Step [II]: $\gamma \subset Z(f)$. If $p\not \in Z(f)$, then $H_{0,\alpha_0}(p)>0$. By continuity, there exist $U$ open of $0$ and $V$ open of $p$ such that $H_{\lambda}(q)>0$ for $q \in V$ and $\lambda \in U$. So $p \notin \gamma$.\\
\\
Step [III]: $\gamma \subset \Gamma$. Take $p \notin \bar{B}_R$ and $U$ a neighborhood of $p$ such that $U \cap \bar{B}_R = \emptyset$. Then, for $q \in U$, and $\lambda>0$ small enough $H_{\lambda}(q)>0$. So $p \notin \gamma$.\\
\\
Step [IV]: $\Gamma \subset \gamma$. Take $p \in \Gamma$ and $\delta>0$. First of all, as $f(x,y)^2\geq0$ and $f(x,y) \not\equiv 0$:
$$
lim_{\lambda \rightarrow 0} \max_{q \in \overline{B_\delta(p)}} H_{\lambda}(q)>0
$$
Now, assume that $p \notin \Gamma_E$. In this case, for $\lambda$ small enough:
$$
H_{\lambda}(p)= \lambda(|p|^2-R^2) \prod_{e \in E}{(|p-p_e|^2-\lambda^2)} - \alpha_0 \lambda^4 < 0
$$
So, for $\lambda$ small enough, there exists $q_{\lambda} \in \overline{B_\delta(p)}$ such that $H_{\lambda}(q_{\lambda})>0$ and $H_{\lambda}(p)<0$. So there exists $q^0_{\lambda} \in \overline{B_\delta(p)}$ such that $H_{\lambda}(q_{\lambda}^0)=0$: $G_{\lambda} \cap \overline{B_\delta(p)} \neq \emptyset$ for all small enough $\lambda$.\\
\\
Now, take $p=p_e$. As $p$ is a regular point of $Z(f)$, it is not isolated and there always exists $q \in B_{\frac{\delta}{2}}(p)$ such that $q \notin \Gamma_E$. Using the first part of for $\frac{\delta}{2}$ and $q$, we conclude that $G_{\lambda} \cap \overline{B_{\delta}(p)} \neq \emptyset$. As $\delta>0$ is arbitrarily small, the result follows. 
\end{proof}
\begin{lemma}
For almost all $\alpha \in [0,1]$, there exist $(\lambda_n)_{n \in \N}$, $\lambda_n>0$ and $\lambda_n \rightarrow 0$ as $n \rightarrow \infty$, such that $G_{\lambda_n}$ is a regular set (i.e.: $\nabla_{x,y}(H_{\lambda_n,\alpha})(p) \neq (0,0)$ for all $p \in G_{\lambda_n}$).
\label{lem:tecr}
\end{lemma}
\begin{proof}
Suppose by absurd that there exists $I \subset [0,1]$ such that:
\begin{itemize}
\item for each $\alpha \in I$, there exists $\rho_{\alpha}>0$ such that: for $0<\lambda<\rho_{\alpha}$, $0$ is not a regular value of $H_{\lambda,\alpha}(x,y)$; 
\item $I$ has positive measure.
\end{itemize}
Now, define $\mathcal{G}_1=\{\alpha \in [0,1]; \frac{1}{2} <\rho_{\alpha}\ \}$ and $\mathcal{G}_n=\{\alpha \in [0,1]; \frac{1}{n+1} <\rho_{\alpha}\leq \frac{1}{n}\}$ for $n>1$. Clearly $\bigcup_{n \in \omega} \mathcal{G}_n=I$. As there are only a countable number of $\mathcal{G}_n$, there must exist an $n_0$ such that the measure of $\mathcal{G}_{n_0}$ is positive. So, fix $\lambda_0<\frac{1}{n_0+1}$: $0$ is a critical value of $H_{\lambda_0,\alpha}(x,y)$ for every $\alpha \in \mathcal{G}_{n_0}$. But $H_{\lambda_0,\alpha}(x,y)=0$ if, and only if $h_{\lambda_0}(x,y)= \alpha \lambda^4_0$. This implies that the set of critical values of $h_{\lambda_0}(x,y)$ has positive measure. But this contradicts Sard theorem: as $h_{\lambda_0}(x,y)$ is an analytic function, the set of critical values must have zero measure. 
\end{proof}
The three lemmas clearly imply the statements of proposition $\ref{prop:123}$.
\section{Proof of Proposition $\ref{prop:4}$} \label{sec:3}
At this section we prove Proposition $\ref{prop:4}$. For now on, we assume hypothesis $[iv]$ and we fix $\alpha \in [0,1]$. We will denote by $\Gamma_{\partial}$ the set of points $p$ that are contained in $\Gamma \cap \partial \bar{B}_R$. By hypothesis, $\Gamma_{\partial}$ is a finite set of points.\\
\\
We start with four local lemmas.
\begin{figure}
  \centering
\includegraphics[scale=1]{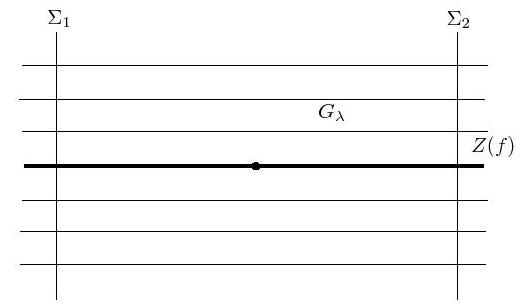}
\includegraphics[scale=1]{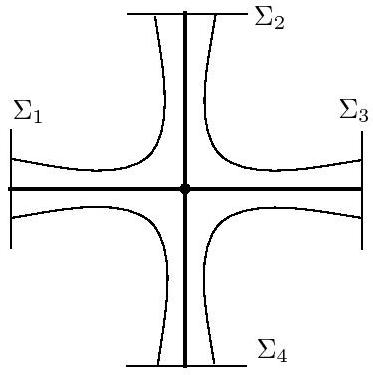}
\caption{\small{At left a representation of $G_{\lambda}$ close to a regular $p$. At right, a representation of $G_{\lambda}$ close to a non-regular $p$ (see lemma $\ref{lem:pontosR}$)}}
\label{fig:0}
\end{figure}
\begin{lemma}
Let $p \in \Gamma$ be a not isolated point of $\Gamma$ such that $p \notin \Gamma_E \cup \Gamma_{\partial}$. Then, there exists a neighborhood $U_p$ of $p$, and $\lambda_p>0$ such that:
\begin{itemize}
\item $Z(f) \cap U$ is the finite union of $\{p\}$ and regular $1$-manifolds $(\gamma_k)_{k \in K}$ with $p$ as a boundary point;
\item $U_p \setminus Z(f)$ is the finite union of connected opens sets $C_l$ for $l \in L$ (in particular $|L|=|K|$);
\item For each $\gamma_k$ there exists a segment $\Sigma_k$ cutting transversally $\gamma_k$;
\item $\Sigma_k$ cuts transversally every $G_{\lambda}$ for $0 \leq \lambda < \lambda_p$;
\item For every pair $(k,l)$ such that $\Sigma_k \cap C_l = \Sigma^l_k \neq \emptyset$, the intersection of $G_{\lambda}$ with $\Sigma^l_k$ is exactly one point for $0<\lambda<\lambda_p$;
\item For each $l \in L$, there exists exactly two $k \in K$ ($k_1$ and $k_2$) such that $\Sigma_k \cap C_l \neq \emptyset$. There exists a transition map $\phi: \Sigma^l_{k_1} \rightarrow \Sigma^l_{k_2}$ such that $q$ and $\phi(q)$ are on the same connected component of $G_{\lambda}$ for every $0<\lambda <  \lambda_p$.
\end{itemize}
see figure $\ref{fig:0}$ for an illustration.
\label{lem:pontosR}
\end{lemma}
\begin{proof}
By the Newton-Puisseux theorem, there exists $U_p$ such that $Z(f) \cap U_p$ is the finite union of $\{p\}$ and regular $1$-manifolds $(\gamma_k)_{k \in K}$ with $p$ as a boundary point. Now, by the implicit function theorem, $\{(x,y,\lambda); H_{\alpha}(x,y,\lambda)\}$ is locally a graph given by a function $\lambda: U_p \rightarrow \R$.\\
\\
Define $C_l$, for $l \in L$ as the connected components of $U_p \setminus Z(f)$. Shrinking $U_p$ if necessary, we can assume that the boundary of each $C_l$ contains exactly two $\gamma_k$. For each $\gamma_k$, chose a point $p_k$ such that $\nabla f(p_k) \neq (0,0)$. Now, there exists $\lambda_p>0$ such that, for $|\lambda|<\lambda_p$:
$$
H_{\lambda,\alpha}(x,y) = f(x,y)^2 - \lambda g(x,y,\lambda)
$$
where $g(x,y,\lambda) = -(x^2+y^2-R^2) \prod_{e \in E}{((x-x_e)^2+(y-y_e)^2-\lambda^2)} + \alpha \lambda^3$ is positive on $U$. Now, take the change of coordinates:
$$
\begin{array}{ccc}
\bar{x} = <\nabla (f(p_k))^{\bot},(x,y)>; &
\bar{y} = \frac{f(x,y)}{\sqrt{g(x,y,\lambda)}}; &
\bar{\lambda}=\lambda;
\end{array}
$$
Close to $p_k$ this change of coordinates is a diffeomorphism that sends $p_k$ to a point $(x_k,0)$. Dropping the bars, we get that:
$$
H_{\lambda,\alpha}(x,y) = U(x,y,\lambda) ( y^2 - \lambda)
$$
where $U(x,y,\lambda)$ is an unit. For some small enough $\delta_k>0$, take $\Sigma^{\ast}_k = \{(x_k,t); t \in ]-\delta,\delta[\}$ and $\Sigma_k$ as the inverse image of $\Sigma^{\ast}_k$ by the local diffeomorphism. It is now clear that $\Sigma_k$ respects all conditions imposed in the enunciate. It rests to show that there exists a transition function.\\
\\
Fixed $l \in L$, take $k_1$ and $k_2$ such that $\Sigma^l_{k_i} \neq \emptyset$. As $\{(x,y,\lambda); H_{\alpha}(x,y,\lambda)\}$ is locally a graph, $G_{\lambda} \cap U_p$ is a level curve. As $\lambda(x,y)=0$ if, and only if, $(x,y) \in Z(f)$, each $G_{\lambda}$ for  $0<\lambda<\lambda_p$ can not cross $Z(f)$.\\
\\
Now, by continuity, shrinking $\lambda_p$ if necessary, the points of $G_{\lambda}$ cutting $\Sigma^l_{k_i}$ are on the same connected component of $G_{\lambda}$. We claim that there can not exist another connected component of $G_{\lambda} \cap C_l$.\\
\\
Indeed, consider the set $S$ contained in between $Z(f)$, $\Sigma^l_{k_i}$ and some connected part of $G_{\lambda_0}$ cutting the two $\Sigma^l_{k_i}$. For $\lambda$ small enough, if $G_{\lambda}$ is disconnected, there exists a connected part of it contained in $S$. As $G_{\lambda}$ is locally a level curve, this component must contain an oval. This implies that there exists a point $q \in S$ such that $\nabla \lambda(q)=(0,0)$. Now, suppose by absurd that there exists $(\lambda_n)_{n \in \N}$ an infinite sequence converging to zero such that $G_{n}$ is disconnected. There would exist an infinite sequence $q_n$ converging to some point of $Z(f)$ such that $\nabla \lambda(q_n)=(0,0)$. As $\lambda(q)$ is analytic, this would imply the existence of a curve inside $S$ where $\lambda(x,y)$ would be zero. But this is an absurd. 
\end{proof}
Now, we treat the isolate points:
\begin{lemma}
If $p \in \Gamma$ is isolated, then there exists $U_p$ neighborhood of $p$ and $\lambda_p>0$ such that $G_{\lambda} \cap U_p$ is connected for all $0<\lambda<\lambda_p$.
\label{lem:pontoI} 
\end{lemma}
\begin{proof}
First, notice that $p$ can not be contained in $\Gamma_E \cup \Gamma_{\partial}$ by definition. After a translation, we may assume that $p=(0,0)$. By the implicit function theorem, there exists $U_p$ neighborhood of $p$ where $\{(x,y,\lambda);H_{\alpha}(x,y,\lambda)=0\}$ is a graph given by a function $\lambda: U_p  \rightarrow \R$. There exists $\lambda_p>0$ such that:
$$
H_{\alpha}(x,y,\lambda) = f(x,y)^2 - \lambda g(x,y,\lambda)
$$
where $g(x,y,\lambda)$ is positive on $U_p$ for $|\lambda|<\lambda_p$. Now, fix $\Sigma=\{(x,0); x \in \R\}$, there exists $n\in \N$ such that:
$$
H_{\alpha}(x,0,\lambda) = x^{2n}\bar{f}(x)^2 - \lambda g(x,0,\lambda)
$$
where $\bar{f}(x)$ is an unit. Then, taking the local diffeomorphism:
$$
\begin{array}{ccc}
\bar{x}= \frac{x \bar{f}(x)^{\frac{2}{2n}}}{g(x,0,\lambda)^{\frac{1}{2n}} }; & \bar{y}=y; & \bar{\lambda} = \lambda;
\end{array}
$$
we get:
$$
H_{\alpha}(x,0,\lambda) = U(x,\lambda)(x^{2n} - \lambda )
$$
where $U(x,\lambda)$ is an unit. Shrinking $U$ if necessary, this implies that, for $0<\lambda<\lambda_p$, $G_{\lambda}\cap \Sigma \cap U$ is equal to two points: one with positive $x$ coordinate, and the other with negative $x$ coordinate. Again shrinking $U$ and $\lambda_p$ if necessary, by lemma $\ref{lem:convergencia}$, we can assume that $G_{\lambda}$ does not pass through the boundary of $U$. As $G_{\lambda}$ must be closed, the two points of $G_{\lambda} \cap \Sigma \cap U$, must be on the same connected component of $G_{\lambda}$.\\
\\
The existence of more connected component is not possible by a similar argument of lemma $\ref{lem:pontosR}$.
\end{proof}
Now we turn ourself to the special points $p \in \Gamma_E$.
%
%
%
%
\begin{figure}
  \centering
\includegraphics[scale=1.2]{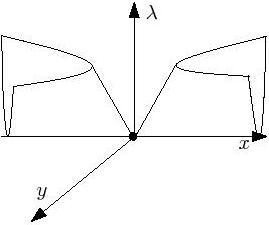}
\includegraphics[scale=1]{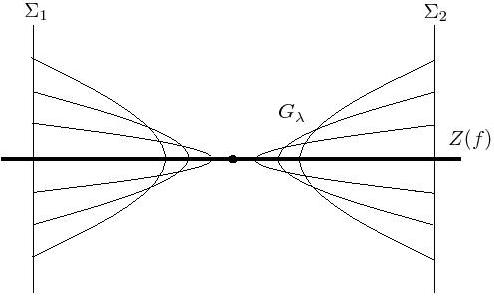}
\caption{\small{At left a representation of the surface $\{H_{\alpha}(x,y,\lambda)=0\}$ close to a point $p_e$. At right, a representation of $G_{\lambda}$ close to $p_e$ (see lemma $\ref{lem:pontosPe}$)}}
\label{fig:1}
\end{figure}

\begin{lemma}
Let $p \in \Gamma_E$. Then, there exists a neighborhood $U_p$ of $p$, and $\lambda_p>0$ such that:
\begin{itemize}
\item $Z(f) \cap U_p$ is a regular curve;
\item There exists a curve $\sigma \subset U_p$ transversal to $Z(f)$ at $p$ such that $U_p \setminus \sigma$ is the union of two connected open sets $C^j$ for $j=1,2$;
\item $C^j \setminus Z(f)$ is the union of two connected sets $C_i^j$ for $i=1,2$;
\item There exists two segments $\Sigma^j \subset C^j $ for $j=1,2$, transversal to $Z(f)$ and to every $G_{\lambda}$ for $0 < \lambda < \lambda_p$. We will denote by $\Sigma_i^j$ the intersection $\Sigma^j \cap C^j_i$;
\item The intersection $\Sigma^j_i \cap G_{\lambda}$, for $0<\lambda<\lambda_p$, is exactly one point;
\item $G_{\lambda} \cap U_p$ is the union of two connected curves for $0 < \lambda < \lambda_p$;
\item There exists two transition maps $\phi^j: \Sigma^j_1 \rightarrow \Sigma^j_2$ such that $q$ and $\phi^j(q)$ are on the same connected component of $G_{\lambda} \cap U_p$ for $0 < \lambda < \lambda_p$.
\end{itemize}
see figure $\ref{fig:1}$ for an illustration.
\label{lem:pontosPe}
\end{lemma}
\begin{proof}
Let us start making a translation and a rotation of $(x,y)$ such that $p$ is send to the origin and $\nabla f (x,y)=(0,c)$, for some $c \neq 0$. At this new coordinates, we have:
$$
H_{\alpha}(x,y,\lambda) = f(x,y)^2  - \lambda(x^2+y^2-\lambda^2) g(x,y,\lambda) -\alpha \lambda^4
$$
where $g(x,y,\lambda)$ is a positive function in a neighborhood of $(0,0,0)$. The change of coordinates:
$$
\begin{array}{ccc}
\bar{x}=x; & \bar{y} = \frac{f(x,y)}{\sqrt{g(x,y,\lambda)}}; & \bar{\lambda}=\lambda;
\end{array}
$$
is clearly a diffeomorphism in a neighborhood of the origin. Dropping the bars, we get:
$$
H_{\alpha}(x,y,\lambda) = y^2  - \lambda(x^2+ C y^2-\lambda^2 + \phi(x,y,\lambda)) -\alpha \lambda^4 \psi(x,y,\lambda)
$$
where $C= \frac{c^2}{g(0,0,0)}\neq0$, $\phi(x,y,\lambda) \in o(x,y,\lambda,3)$, and $\psi(0,0,0) \neq 0$. Taking the expansions:
$$
\begin{array}{c}
\phi(x,y,\lambda)=y^2\phi_2(x,y,\lambda)+y\phi_1(x,\lambda)+\phi_0(x,\lambda)\\
\psi(x,y,\lambda)=y^2\psi_2(x,y,\lambda)+y\psi_1(x,\lambda)+\psi_0(x,\lambda)
\end{array}
$$
where $\phi_2(x,y,\lambda) \in o(x,y,\lambda,1)$, $\phi_1(x,\lambda) \in o(x,y,\lambda,2)$, $\phi_0(x,\lambda) \in o(x,y,\lambda,3)$ and $\psi_i \in o(x,y,\lambda,0)$ for all $i=0,1,2$. We get:
$$
\begin{array}{c}
H_{\alpha}(x,y,\lambda) = y^2\left(1-\lambda(C+\phi_2(x,y,\lambda)+\alpha \lambda^3 \psi_2(x,y,\lambda))\right) +\\
+ y\left(\lambda(\phi_1(x,\lambda)+\alpha \lambda^3 \psi_1(x,\lambda))\right) - \lambda \left(x^2-\lambda^2+\phi_0(x,\lambda)+\lambda^3 \psi_0(x,\lambda)\right)
\end{array}
$$
In particular, the expression multiplying $y^2$ is non-zero at the origin (it is equal to $1$). So, we can apply the Weierstrass preparation theorem to get:
$$
H_{\alpha}(x,y,\lambda) = \mu(x,y,\lambda)(y^2+a(x,\lambda)y + b(x,\lambda))
$$
where $\mu(0,0,0) \neq 0 $. Taking the expansion:
$$
\mu(x,y,\lambda)=y^2\mu_2(x,y,\lambda)+y\mu_1(x,\lambda)+\mu_0(x,\lambda)
$$
we get the system:
$$
\begin{array}{c}
\mu_0(x,\lambda)b(x,\lambda) = \lambda\left(x^2-\lambda^2+\phi_0(x,\lambda)+\lambda^3 \psi_0(x,\lambda)\right)\\
\mu_0(x,\lambda)a(x,\lambda) + \mu_1(x,\lambda)b(x,\lambda) = 
\lambda(\phi_1(x,\lambda)+\alpha \lambda^3 \psi_1(x,\lambda))\\
\mu_0(x,\lambda) + \mu_1(x,\lambda) a(x,\lambda) + \mu_2(x,0,\lambda)b(x,\lambda)=\left(1-\lambda(C+\phi_2(x,0,\lambda)+\alpha \lambda^3 \psi_2(x,0,\lambda))\right)
\end{array}
$$
As $\mu_0(0,0) \neq 0$, from the system we get that $\mu_0(0,0)=1$. So:
$$
\begin{array}{c}
b(x,\lambda)=\lambda\left(x^2-\lambda^2+B(x,\lambda) \right)\\
a(x,\lambda)=\lambda A(x,\lambda)
\end{array}
$$
where $B(x,\lambda) \in o(x,\lambda,3)$ and $A(x,\lambda) \in o(x,\lambda,2)$. So, the surface $\{H_{\alpha}(x,y,\lambda) =0\}$ is locally given by:
$$
y^2+\lambda A(x,\lambda)y + \lambda \left(x^2-\lambda^2+B(x,\lambda) \right) = 0 
$$ 
Solving the equation we get:
$$
y= -\lambda \frac{A(x,\lambda)}{2} \pm \sqrt{\lambda^2 \frac{A^2(x,\lambda)}{4}-\lambda\left(x^2-\lambda^2+B(x,\lambda) \right)}
$$
which can be re-written as:
$$
y= \Phi(x,\lambda) \pm \sqrt{\lambda( \lambda^2 - x^2+ \Psi(x,\lambda) )}
$$
where $\Phi(x,\lambda) \in o(x,\lambda,3)$ and $\Psi(x,\lambda) \in o(x,\lambda,3)$. Now $\{\lambda^2 - x^2+ \Psi(x,\lambda)=0\}$ are the points where there exists only one solution of $y$. As the Hessian of this function is not degenerated, we conclude that $\{\lambda^2 - x^2+ \Psi(x,\lambda)=0\}$ is given, after a diffeomorfism tangent to the identity, by $\{(\lambda - x)(\lambda+x) =0\}$.\\
\\
So, in a small neighborhood $V$ of the origin, there exists two continuous curves $(C_1(\lambda),\lambda)$ and $(C_2(\lambda),\lambda)$ such that:
\begin{itemize}
\item $C_i(\lambda) \rightarrow 0$ for $i=1,2$, when $\lambda \rightarrow 0$;
\item $C_1(\lambda)>0 $ and $C_2(\lambda)<0$ for $\lambda>0$;
\item for $x>C_1(\lambda)$ or $x<C_2(\lambda)$, there exists two $y$ coordinates such that $(x,y,\lambda) \in H_{\alpha}^{-1}(0)$;
\item for $x=C_1(\lambda)$ or $x=C_2(\lambda)$, there exists one $y$ coordinate such that $(x,y,\lambda) \in H_{\alpha}^{-1}(0)$;
\item for $C_2(\lambda)<x<C_1(\lambda)$, there does not exist a $y$ coordinate such that $(x,y,\lambda) \in H_{\alpha}^{-1}(0)$.
\end{itemize}
The existence of $C_1(\lambda)$ and $C_2(\lambda)$ grantees that $G_{\lambda}$ is locally the union of two disconnected curves.\\
\\
By a similar reasoning it is easy to prove the existence of two continuous curves $(c_1(\lambda),\lambda)$ and $(c_2(\lambda),\lambda)$ such that: 
\begin{itemize}
\item $c_i(\lambda) \rightarrow 0 $ for $i=1,2$, when $\lambda \rightarrow 0$;
\item $c_1(\lambda)\geq C_1(\lambda) >0$ and $c_2(\lambda) \leq C_2(\lambda)<0$
\item for $x>c_1(\lambda)$ or $x<c_2(\lambda)$, there exists a positive and a negative $y$ coordinate such that $(x,y,\lambda) \in H_{\alpha}^{-1}(0)$.
\end{itemize}
At this new coordinates, the lemma is clear.
%
%
\end{proof}
We now study the points $p \in \Gamma_{\partial}$:
\begin{figure}
  \centering
\includegraphics[scale=1.2]{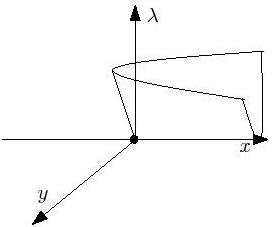}
\includegraphics[scale=1]{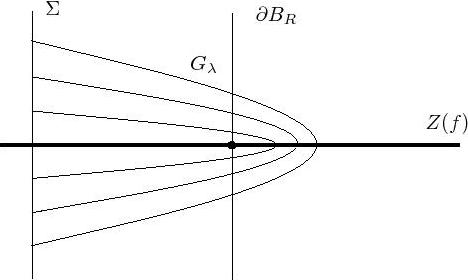}
\caption{\small{At left a representation of the surface $\{H_{\alpha}(x,y,\lambda)=0\}$ close to a point $p_a$. At right, a representation of $G_{\lambda}$ close to $p_a$ (see lemma $\ref{lem:pontosPb}$)}}
\label{fig:2}
\end{figure}
\begin{lemma}
Let $p \in \Gamma_{\partial}$. Then, there exists a neighborhood $U_p$ of $p$, and $\lambda_p>0$ such that:
\begin{itemize}
\item $Z(f) \cap U_p$ is a regular curve and $U_p \setminus Z(f)$ is the union of two connected components $C_i$ for $i=1,2$;
\item There exists a segment $\Sigma \subset B_R$, transversal to $Z(f)$ and to every $G_{\lambda}$ for $0 < \lambda < \lambda_p$. We will denote by $\Sigma_i$ the intersection $\Sigma \cap C_i$;
\item The intersection $\Sigma_i \cap G_{\lambda}$ is exactly one point for $0<\lambda<\lambda_p$;
\item $G_{\lambda} \cap U_p$ is a connected curve for $0 < \lambda < \lambda_p$;
\item There exists a transition map $\phi: \Sigma_1 \rightarrow \Sigma_2$ such that $q$ and $\phi_j(q)$ are on the same connected component of $G_{\lambda} \cap U_p$ for $0 < \lambda < \lambda_p$.
\end{itemize}
see figure $\ref{fig:2}$ for an illustration.
\label{lem:pontosPb}
\end{lemma}
\begin{proof}
For shortness, we will omit all calculus that are analogous with the last proof. Take a translation that sends $p$ to the origin and a change of coordinates such that, locally, $\partial B_R$ is now given by $\{x=0\}$:
$$
H_{\alpha}(x,y,\lambda) = f(x,y)^2  + \lambda x g(x,y,\lambda) -\alpha \lambda^4
$$
where $g(x,y,\lambda)$ is a positive function in a neighborhood of $(0,0,0)$. Take the change of coordinates:
$$
\begin{array}{ccc}
\bar{x} = x; & \bar{y}=\frac{f(x,y)}{\sqrt{g(x,y,\lambda)}}; & \bar{\lambda}=\lambda;
\end{array}
$$
It is a local diffeomorphism because $\nabla f(p)$ is transversal to $\partial B_r$. Dropping the bars we get:
$$
y^2 + \lambda x - \alpha \lambda^4 \phi(x,y,\lambda)
$$
where $\phi(0,0,0) \neq 0$. Re-writing, we get:
$$
y^2(1-\alpha \lambda^4 \phi_2(x,y,\lambda)) - y \alpha \lambda^4 \phi_1(x,\lambda) + \lambda(x-\alpha \lambda^3 \phi_0(x,\lambda))
$$
By the Weierstrass preparation theorem, we get:
$$
H_{\alpha}(x,y,\lambda) = \mu(x,y,\lambda)(y^2+a(x,\lambda)y + b(x,\lambda))
$$
and by an analogous argument as of the other proof, we get:
$$
\begin{array}{c}
b(x,\lambda)=\lambda (x + B(x,\lambda))\\
a(x,\lambda) = \lambda A(x,\lambda)
\end{array}
$$
where $B(x,\lambda) \in o(x,\lambda,2)$ and $A(x,\lambda) \in o(x,\lambda,1)$. Solving in function of $y$ we get:
$$
y= \Phi(x,\lambda) \pm \sqrt{\lambda( -x + \Psi(x,\lambda) )}
$$
where $\Phi(x,\lambda) \in o(x,\lambda,2)$ and $\Psi(x,\lambda) \in o(x,\lambda,2)$. So, by the implicit function theorem, for $\lambda$ small enough, there exists a continuous curve $(C(\lambda),\lambda)$ such that:
\begin{itemize}
\item $C(\lambda) \rightarrow 0 $ when $\lambda \rightarrow 0$;
\item for $x<C(\lambda)$, there exists two $y$ coordinates such that $(x,y,\lambda) \in H_{\alpha}^{-1}(0)$;
\item for $x=C(\lambda)$, there exists one $y$ coordinates such that $(x,y,\lambda) \in H_{\alpha}^{-1}(0)$;
\item for $x>C(\lambda)$, there does not exist a $y$ coordinates such that $(x,y,\lambda) \in H_{\alpha}^{-1}(0)$.
\end{itemize}
More than this, again by the implicit function theorem, there exists a curve $c(\lambda)$ such that:
\begin{itemize}
\item for $x<c(\lambda)$, there exists a positive and a negative $y$ coordinate such that $(x,y,\lambda) \in H_{\alpha}^{-1}(0)$.
\end{itemize}
On this new coordinates, the lemma is clear.
\end{proof}
After these four local lemmas, we are ready to prove proposition $\ref{prop:4}$:
\begin{proof}\textit{(Proposition \ref{prop:4})}
For each point $p \in \Gamma$, we have a local description of $G_{\lambda}$ at an open set $U_p$ and $0 \leq \lambda< \lambda_p$ by lemmas $\ref{lem:pontosR}$, $\ref{lem:pontoI}$, $\ref{lem:pontosPe}$ and $\ref{lem:pontosPb}$. As $\Gamma$ is compact, there exists a finite number of points $p_k$ for $k \in K$ such that $\Gamma \subset \cup_{k \in K}U_{p_k}$. Take $\lambda_0= min\{\lambda_{p_k}; k \in K\}$. By lemma $\ref{lem:convergencia}$, we can assume that $\lambda_0$ is small enough so $G_{\lambda} \subset \cup_{k \in K}U_{p_k}$.\\
\\
Let us assume that $\Gamma$ is connected. If $\Gamma$ is only one point then the proposition is clear by lemma $\ref{lem:pontoI}$. If not, 
let us first do an `intuitive' argument to clear the idea of the proof. Recall the definition of support of $(f,R)$ on section $\ref{sec:1}$ and take $C_i$ an internal component. If there was no point $p \in \Gamma_E$ on the boundary of $C_i$, then, by lemma $\ref{lem:pontosR}$, $\Gamma_{\lambda} \cap C_i$ would be connected and would not cross $\partial C_i$. But this is not the case, as all $C_i$ internal must have a point $p \in \Gamma_E$ on its boundary. Such a point, by lemma $\ref{lem:pontosPe}$ `intuitively' glue two ovals of different $C_i$'s. So, if we glue all ovals contained in each $C_i$ in a `good manner', we will get a single connected component.\\
\\
Now, let us formalize the idea. First, take $p_e \in \Gamma_E$ for $e=<i,j>$. For every $0<\lambda<\lambda_0$, by lemma $\ref{lem:pontosPe}$, $G_{\lambda}$ is the union of two locally disconnected curves $g_1$  $g_2$. We claim that this curves got to be part of the same connected component of $G_{\lambda}$. Indeed, each one of these curves must be part of a connected oval of $G_{\lambda}$, and thus have to be closed. But they are locally disconnected by $\Gamma$: $g_1$ and $g_2$ has a part contained in $C_i$ and another in $C_j$. Now, $G_{\lambda}$ can only cross $\Gamma$ in a neighborhood of points $\Gamma_E \cup \Gamma_{\partial}$. So, for the part of $g_1$ that is contained in $C_j$ to reach the $C_i$ part, it would have to pass through a finite number of $C_{k_n}$'s and the $p_{e_n}$ for $e_n=<k_n,k_{n+1}>$ where $k_1=j$ and $k_N=i$. But we have chosen $T$ to be a tree, and thus have no cycle, the only possibility is to have $k_{N-1}=j$. As $p_e$ is unique for each $E$, this can only happen if $g_1$ and $g_2$ are `globally connected': they are on the same connected component of $G_{\lambda}$.\\
\\
Now, take $C_i$ an interior component. We claim that $G_{\lambda} \cap C_i$ is contained in the same connected component of $G_{\lambda}$. Indeed, the locally disconnected curves close to $p_e \in (\Gamma_E \cap \partial C_i)$ are globally connected. Using lemma $\ref{lem:pontosR}$, we conclude that there exists transitions between two neighbors $p_e \in (\Gamma_E \cap \partial C_i)$, and thus they are on the same connected component of $G_{\lambda}$. By transitivity, we conclude the claim.\\
\\
\begin{figure}
  \centering
  \includegraphics[scale=2]{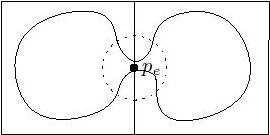}
  \caption{\small{Illustration on how two ovals are glued near a point $p_e \in \Gamma_E$.}}
 \end{figure}
It is now easy to see that this argument is true for $C_0$, if we first notice that lemma $\ref{lem:pontosPb}$ connects the ovals of two disconnected components of $C_0$. Also, it is clear that if $\{i,j\} \in E$, then $\Gamma \cap (C_i \cup C_j)$ is contained in the same connected component of $G_{\lambda}$. By transitivity, this is true for any $i,j$ such that there exists a walk from $i$ to $j$. As $T$ is a tree, and thus connected, all $G_{\lambda}\cap C_i$ are contained in the same connected component of $G_{\lambda}$ and $G_{\lambda}$ is connected.\\
\\
If $\Gamma$ is the union of $N$ connected components, then we can assume that the neighborhoods $U_{p_i}$ are small enough so that $\cup U_{p_i}$ is the union of $N$ connected open sets $V_{n}$, each one containing a different connected part of $\Gamma$. Restraining the tree of adjacency $T$ and the support of $(f,R)$ to each one of this connected open sets $V_{n}$ and using the first part of the proof, we conclude that $G_{\lambda} \cap V_{n}$ is connected for each $n \leq N$. So, $G_{\lambda}$ is an union of $N$ connected components, each one of them close to a different connected component of $\Gamma$. 
\end{proof}
\section*{Aknowledgments}
I would like to thank professor Daniel Panazzolo for the usefull talks, all suggestions and revision. This work was supported by the Universit\'e de Haute Alsace.
\end{document}